\documentclass[
11pt,%
final,%
tightenlines,%
notitlepage,%
twoside,%
onecolumn,%
nobibnotes,%
nofootinbib,%
superscriptaddress,%
noshowpacs,%
noshowydk,%
showdoi,%
centertags,%
maik]%
{revtex4n}

\usepackage{listings}

\input jpack_e
\setcaptionmargin{5mm}
\onelinecaptionstrue

\columntovsizetrue
\allowdisplaybreaks

\begin{document}

  \newcommand{\rtxt}[1]{{\textcolor{red} {#1}}}
\newcommand{\btxt}[1]{{\textcolor{blue} #1}}
\newcommand{\AT}[0]{\textcolor{black}{AT }}
\newcommand{\GC}[0]{\textcolor{black}{GC }}
\newcommand{\R}{\mathbb{R}}
\newcommand{\K}{\mathbb{K}}
\newcommand{\N}{\mathbb{N}}

\newtheorem{teo}{Theorem}
\newtheorem*{teon}{Theorem}
\newtheorem{lem}{Lemma}
\newtheorem*{lemn}{Lemma}
\newtheorem{prp}{Proposition}
\newtheorem*{prpn}{Proposition}
\newtheorem{ass}{Assertion}
\newtheorem*{assn}{Assertion}
\newtheorem{assum}{Assumption}
\newtheorem*{assumn}{Assumption}
\newtheorem{stat}{Statement}
\newtheorem*{statn}{Statement}
\newtheorem{cor}{Corollary}
\newtheorem*{corn}{Corollary}
\newtheorem{hyp}{Hypothesis}
\newtheorem*{hypn}{Hypothesis}
\newtheorem{con}{Conjecture}
\newtheorem*{conn}{Conjecture}
\newtheorem{dfn}{Definition}
\newtheorem*{dfnn}{Definition}
\newtheorem{problem}{Problem}
\newtheorem*{problemn}{Problem}
\newtheorem{notat}{Notation}
\newtheorem*{notatn}{Notation}
\newtheorem{quest}{Question}
\newtheorem*{questn}{Question}

\theorembodyfont{\rm}
\newtheorem{rem}{Remark}
\newtheorem*{remn}{Remark}
\newtheorem{exa}{Example}
\newtheorem*{exan}{Example}
\newtheorem{cas}{Case}
\newtheorem*{casn}{Case}
\newtheorem{claim}{Claim}
\newtheorem*{claimn}{Claim}
\newtheorem{com}{Comment}
\newtheorem*{comn}{Comment}

\theoremheaderfont{\it}
\theorembodyfont{\rm}

\newtheorem{proof}{Proof}
\newtheorem*{proofn}{Proof}

\selectlanguage{english}
\Rubrika{\relax}
\CRubrika{\relax}
\SubRubrika{\relax}
\CSubRubrika{\relax}
%

\def\JournalNumber{2}
\def\JournalVolume{23}
%
%
%
\nameVolumeRus{}
\CnameVolumeRus{}
\nameIssueRus{\No}
\CnameIssueRus{}
\namePartRus{}
\namePagesRus{}
\nameYearShortRus{}
\JournalNameRus{}
\TranslitJournalNameRus{}
\JournalName{Regular and Chaotic Dynamics}
\JournalISSNCode{1560-3547}
\IssuePrice{}
\TransYearOfIssue{2018}
\TransCopyrightYear{2016}%
\OrigYearOfIssue{}
\OrigCopyrightYear{2018}%
\OrigIssueNo{\JournalNumber}
\OrigVolumeNo{\JournalVolume}
\TransVolumeNo{\JournalVolume}
\TransIssueNo{\JournalNumber}
\TransPartNo{}
\SHORTjournalPREFIX{RCD} 
\LONGjournalPREFIX{RegDyn} 
\BatFileName{call make_ps.bat} 
\BatSwitch{3} 
\IssueName{}
\SupplementNumber{}
\PublicationSerialNumberInYear{0}
\PublicationSerialNumberInVolume{0}
\ConditionalIssueDate{"year","month","day","name","type"}
\PagePrefix{}
\JournalISSNonlineCode{}
\JournalISSNCodeRus{}
\JournalISSNonlineCodeRus{}
\VolumeName{}
\IssnoName{none}
\PartnoName{}
\FpageNamepp{}
\FpageNnamep{}
\FpagePrefix{}
\LpageNnamepp{}
\LpageNamep{}
\LpagePrefix{}
\VolumePageNumbering{}
\JournalPubID{}
\FirstJournalPageNumber{}
\LastJournalPageNumber{}
\makeatletter
\def\MAIKlogo{Pleiades~Publishing,~Ltd.}
\def\maikpraefix{10.0000/S}
\edef\@ContentsHeadLineB{Simultaneous English language translation of the journal is available from \noexpand\MAIKlogo}
\def\Distributed{Distributed worldwide by Springer. }
\def\ArticlePages#1{\relax}
\@ifxundefined\CONT@sw{\@booleantrue\CONT@sw}{}%
\@booleantrue\showPACS@sw%
\@booleantrue\showKEYS@sw %
\@booleantrue\noOrigJournalVersion@sw
\@booleantrue\noOrigVolumeNo@sw
\@booleanfalse\noTransVolumeNo@sw
\makeatother
\input maikdoi %

\beginpaper


\input engnames
\titlerunning{P\'olya Counting in Periodic DNA Chains}
\authorrunning{Hillebrand et al.}
\toctitle{P\'olya Counting in DNA Chains}
\tocauthor{M.\, Hillebrand}
\title{Distribution of Base Pair Alternations in a Periodic DNA Chain: Application of P\'olya Counting to a Physical System}
\firstaffiliation{
}%
\articleinenglish 
\PublishedInRussianNo
\author{\firstname{Malcolm }~\surname{Hillebrand}}%
\email[E-mail: ]{malcolm.hillebrand@gmail.com}
\affiliation{
Department of Mathematics and Applied Mathematics, University of Cape Town,\\
Rondebosch, Cape Town 7701, South Africa}%
\author{\firstname{Guy}~\surname{Paterson-Jones}}%
\affiliation{
Department of Mathematics and Applied Mathematics, University of Cape Town,\\
Rondebosch, Cape Town 7701, South Africa}%
\author{\firstname{George}~\surname{Kalosakas}}%
\affiliation{
Department of Materials Science, University of Patras,
Rio GR-26504, Greece}%
\author{\firstname{Charalampos}~\surname{Skokos}}%
\affiliation{
Department of Mathematics and Applied Mathematics, University of Cape Town,\\
Rondebosch, Cape Town 7701, South Africa}%
\begin{abstract}
In modeling DNA chains, the number of alternations between Adenine-Thymine (AT) and Guanine-Cytosine (GC) base pairs can be considered as a measure of the heterogeneity of the chain, which in turn could affect its dynamics.  {A probability distribution function of the number of these alternations is derived for circular or periodic DNA}. Since there are several symmetries to account for in the periodic chain, necklace counting methods are used. In particular, P\'olya's Enumeration Theorem is extended for the case of a group action that preserves partitioned necklaces. This, along with the treatment of generating functions as formal power series, allows for the direct calculation of the number of possible necklaces with a given number of AT base pairs, GC base pairs and  alternations. The theoretically obtained probability distribution functions of the number of alternations are accurately reproduced by Monte Carlo simulations and fitted by Gaussians. The effect of the number of base pairs on the characteristics of these distributions is also discussed, as well as the effect of the ratios of the numbers of AT and GC base pairs.
\end{abstract}
\received{October 13 2017}
\accepted{December 11, 2017}%
\maketitle

\textmakefnmark{0}{)}%

\section{Introduction}
\label{sec:Intro}

{Single circular DNA molecules are abundant in nature. The whole genome in a typical bacterium
is usually contained in a closed DNA molecule, while in eucaryotes the organelle DNA, inside the
mitochondria and chloroplasts, is also found in the same form \cite{r1,r2}. Also plasmids, either
naturally found in bacteria, or used as vectors in gene cloning, are smaller circular DNA segments.
Apart from these cases, } in considering the dynamics and other properties of DNA chains, it is often useful to model the chain using periodic boundary conditions in order to avoid finite size or edge effects. For example, periodic boundary conditions have been used to study denaturation bubbles and the melting behavior of DNA \cite{DPB93,aresPRL,theodPRE,boianNAR9,faloPRE,zoliJPCM}, probability distributions of thermal openings in the double strand \cite{NL07,JCP}, bubble opening profiles in promoter regions which regulate gene transcription \cite{NAR,EPL,choiBJ,boianPLOSCB,boianNAR10,PB,huangJBE}, binding sites of DNA-associated proteins \cite{nowakPLOSCB,faloPLOSCB}, various dynamical and nonlinear properties of DNA \cite{yakush,peyrNonlin,farago,NL04,CPL,zoliJTB}, as well as charge transport in DNA \cite{hennigEPJB,PRE5,PRE11,kofanePhysA,velardeEPJB}.

\begin{figure}[h!]
  \centering
  \includegraphics{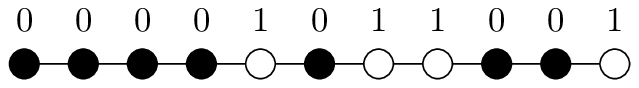}
  \label{lol}
  \caption{An example of a DNA chain. GC base pairs are represented by black beads and the number 0, while AT base pairs are represented by white beads and the number 1. In the DNA necklace corresponding to this chain, the AT base pair at the far right neighbors the GC base pair at the far left.}
\end{figure}

\begin{figure}[h!]
  \centering
  \includegraphics{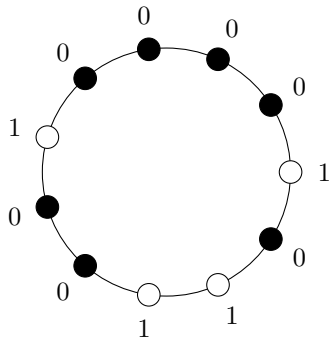}
  \caption{The DNA necklace corresponding to the chain of Fig.~\ref{lol}. This necklace has $\alpha=6$ alternations.}
  \label{circlydna}
\end{figure}

{A DNA chain consists of a series of base pairs, where each base pair is either Adenine-Thymine (AT) or Guanine-Cytosine (GC). Currently, we are investigating the influence of different factors on the chaoticity of periodic DNA chains \cite{Skokos_prep}. One of the examined quantities is the number of base pair alternations, which can be considered as a quantifier of the system's heterogeneity. In this work we focus on the rigorous mathematical treatment of alternation counting in periodic DNA sequences.}
{To study periodic DNA, we will consider the DNA necklace associated to a DNA chain}, where the first and the last base pairs in the chain will become neighbors. This periodicity presents some modeling challenges - if one considers two distinct chains of DNA, it may still be the case that their corresponding necklaces are the same, as one may be merely a rotation or reflection of the other. Such symmetries need to be addressed if any conclusions are to be made about the structure and the dynamics of DNA necklaces. In particular, we are concerned with the number $\alpha$ of base pair alternations in the necklace, where an alternation is defined to be a point at which an AT base pair neighbors a GC base pair or vice versa. Consider, for instance, the DNA chain shown in Fig.~\ref{lol}.
Representing a GC base pair (black bead) with a 0 and an AT base pair (white bead) with a 1, the chain can be written in the form $(1)000\bar0\bar1\bar01\bar10\bar0\bar1(0)$. Here, we have given the leftmost base pair at each alternation point an overbar, and used brackets to denote the fact that in the corresponding DNA necklace the first and last base pairs are neighbors. This necklace is illustrated in Fig.~\ref{circlydna}, and counting the number of overbars we see that there are $\alpha = 6$ alternations.

{It is worth noting that a base pair alternation corresponds to the appearance of the particular sequences (often referred to as ``words'') 01 or 10 in a DNA chain. Word occurrence probabilities have already been studied in the literature (see e.g.~\cite{Kolpakov_etal_2003_NuclAcidsRes_31_3672, Li_ComputChem_1997_21_257, Regnier_2000_104_259, Robin_Daudin_1999_JApplProb_36_179,Robin_Schbath_2001_JCompBio_8_349, Schbath_1995_ESAIM_ProbStat_1_1, Schbath_etal_1995_JCompBio_2_417} and references therein), with emphasis on the appearance of patterns with unexpectedly high or low frequencies, as well as on repeating sequences. However these studies concern the case of linear DNA segments, or in other words DNA chains with fixed boundary conditions. The periodic boundary conditions we consider in our study make the problem of counting alternations (or more generally the appearance of specific words) in circular DNA segments much more complicated than in the case of linear DNA segments due to the appearance of additional symmetries in the DNA structures imposed by rotations and/or reflections.}

Each base pair in a DNA necklace can contribute at most $2$ alternations, depending on which neighbors it differs from. Supposing that the number of AT and GC base pairs in the necklace is given by $N_{AT}$ and $N_{GC}$ respectively, this yields the restriction $0 \leq \alpha \leq \textrm{min}\{2N_{AT}, 2N_{GC}\}$.
We note that in the extreme case of a homogeneous chain composed of base pairs of the same kind $\alpha=0$, while if both types of base pairs are present in the DNA chain the smallest possible value of alternations is $\alpha=2$. The later corresponds to a chain having all AT (and consequently GC) base pairs grouped together. Furthermore, if we traverse the necklace pair by pair until we end up where we started, we must necessarily switch between AT and GC base pairs an even number of times. Thus $\alpha = 2M$ for some $M \in \N$.

Now the natural question is: what is the probability that a random DNA necklace with a specified number of AT and GC base pairs, $N_{AT}$ and $N_{GC}$ respectively, has a specified number of alternations $\alpha$? Or in other words, how many possible combinations of such base pairs are there that yield $\alpha$ alternations once the cyclic and reflective symmetries are taken into account? In what follows we answer these questions and provide an algorithm for computing the number of distinct DNA necklaces satisfying these constraints.

The paper is organized in the following way: In Sect.~\ref{sec:theory}, the mathematical background is laid out, leading into a P\'olya Enumeration Theorem for bipartite sets. In Sect.~\ref{sec:algorithm} an explicit algorithm for calculating the number of distinct DNA necklaces with given values of $\alpha$, $N_{AT}$ and $N_{GC}$ is described, while in Sect.~\ref{sec:NumRes} we compare the theoretical results to those obtained from Monte-Carlo simulations and investigate the effect of the $N_{AT}$ and $N_{GC}$ values on the characteristics of the probability distribution function (pdf) of $\alpha$. Finally, in Sect.~\ref{sec:con} we summarize our results, while in the Appendix we provide a Python computer code implementing the algorithm of Sect.~\ref{sec:algorithm}.

\section{Theoretical Treatment}
\label{sec:theory}

\begin{figure}[h!]
  \centering
  \includegraphics{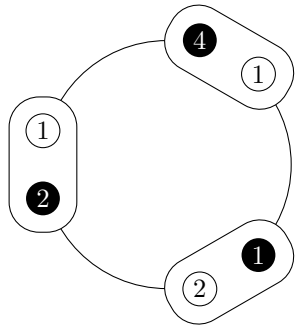}
  \caption{The necklace of containers corresponding to the DNA necklace of  Fig.~\ref{circlydna}. The numbers in each container represent the number of consecutive black or white beads in that segment of the necklace.}
  \label{containerdna}
\end{figure}

Our problem can be neatly related to the combinatorics of necklaces. Effectively, we are interested in the number of distinct necklaces with $N = N_{AT}+N_{GC}$ beads, where $N_{AT}$ of the beads are white, $N_{GC}$ of the beads are black, and there are $\alpha$ alternations between the colors. We consider necklaces to be the same if they can be reflected or rotated into one another, and beads of the same color are treated as indistinguishable. Because of this, we can equivalently think of a necklace with $\alpha$ alternations as a necklace of $\alpha$ \emph{containers}, where each container carries some number of black or white beads of the same color, and adjacent containers have different colors. This idea is illustrated in Fig.~\ref{containerdna}.

We will refer to containers carrying black beads as black containers, and similarly for white containers. Counting the number of distinct necklaces with the given constraints can thus be reformulated as the problem of assigning numbers of beads to $\alpha$ containers, such that the total of the numbers in the black and white containers is equal to $N_{GC}$ and $N_{AT}$ respectively. Two such assignments will be considered equivalent if the containers can be rotated or reflected into one another in such a way as to preserve both the colors and numbers of beads they contain.

Enumerating such assignments is simpler than enumerating necklaces, as we have one less constraint - the number of alternations is now implicit in the formulation of the problem. To perform this enumeration we will require some tools from P\'olya counting theory - in particular, we will need a version of the P\'olya Enumeration Theorem for sets partitioned into two parts, which we will refer to as \emph{bipartite} sets. For completeness' sake, we present this material below.

\subsection{Group Actions}
\label{sec:pergr}

Let $A$ be a set. Then we define the \emph{symmetric group on $A$} to be the set of permutations of $A$:
\begin{equation}
	S_A = \{ \varphi: A \to A \mid \varphi \text{ is a bijection} \}.
\label{eq:Sn}
\end{equation}
A \emph{cycle} is a permutation $\varphi \in S_A$ such that there exist distinct elements $\{ x_1, x_2, \ldots, x_k \} \in A$ and:
\begin{equation}
	\varphi(x) = \begin{cases}
	                x_{i+1} & \text{if $x=x_i$ for some $1 \leq i < k$}\\
	                x_1 & \text{if $x=x_k$}\\
	                x & \text{otherwise.}
                 \end{cases}
\end{equation}
We denote such a cycle suggestively as $(x_1\ x_2\ \ldots\ x_k)$, and say that $\varphi \in S_A$ is a \emph{$k$-cycle} if $\varphi = (x_1\ x_2\ \ldots\ x_k)$ for some $x_i \in S_A$. Two cycles $(x_1\ x_2\ \ldots\ x_k)$ and $(y_1\ y_2\ \ldots\ y_l)$ are said to be \emph{disjoint} if the sets $\{x_1, x_2, \ldots, x_k\}$ and $\{y_1, y_2, \ldots, y_l\}$ are disjoint.

If $A$ is a finite set, every element of $S_A$ can be written as a composition of cycles; in general, however, this cannot be done uniquely. On the other hand, we have the following fundamental structure theorem for elements of finite symmetric groups (see for example \cite{herstein}):
\begin{teon}[Cycle Decomposition Theorem]
	If $A$ is a finite set, then every element $\varphi \in S_A$ can be written as a product of pairwise disjoint cycles, unique up to order of the cycles:
	\begin{equation*}
	    \varphi = (x_{11}\ x_{12}\ \ldots\ x_{1k_1})\ \cdots\ (x_{n1}\ x_{n2}\ \ldots\ x_{nk_n}).
    \end{equation*}
\end{teon}

Given a group $G$ and a set $A$, a \emph{group action} of $G$ on $A$ is a homomorphism $\Gamma_G: G \to S_A$. In other words, elements of $G$ are identified with permutations of $A$ in a manner that preserves the group structure. To simplify the notation, we will write $gx$ instead of $\Gamma_G(g)(x)$ for the action of $g \in G$ on some $x \in A$.

The \emph{orbit} of an element $x \in A$ under the group action $\Gamma_G$ is defined to be the set $\text{Orb}_x = \{gx\ \mid\ g \in G\}$, and its \emph{stabilizer} is given by the subgroup $\text{Stab}_x = \{g \in G\ \mid\ gx = x\}$. Given some $g \in G$, we denote its set of fixed points by $\text{Fix}_g = \{ x \in A\ \mid\ gx = x\}$.

\subsection{P\'olya's Counting Theory}
\label{sec:PCT}

One can often rephrase counting problems in terms of computing the number of distinct orbits of some group action. P\'olya's counting theory can be thought of as a tool for making these computations systematic and expedient. A  fundamental lemma on which this theory is built is the following \cite{burnside}:
\begin{lem}[Burnside's Lemma]
	The number of distinct orbits in a group action of a finite group $G$ on $A$ is given by the average number of fixed points of elements of $G$:
	\begin{equation}
	    \text{\#Orbits} = \frac1{|G|}\sum_{g \in G}|\text{Fix}_g|.
    \label{eq:orb_def}
    \end{equation}
\end{lem}

A basic problem in combinatorics is the  following. Suppose one has a finite set of objects $A$, and one wishes to color them with colors from another set $\Omega$. How many distinct ways are there of coloring the objects up to some kind of symmetry?  This can be recast in the language of group actions. The set of possible colorings is given by $\Omega^{A} = \{\varphi: A \to \Omega\ |\ \varphi \text{ a function} \}$, and the symmetry is given by a group action $\Gamma _G$ on $A$. This group action passes naturally to a group action $\tilde \Gamma _G$ on $\Omega^{A}$, defined by $g\varphi: x \mapsto \varphi(gx)$.

The question now reduces to counting the number of distinct orbits of this latter action. In this simplified case, Burnside's lemma is often sufficient to answer the question. We can generalize this problem slightly, however. Suppose that each color has an associated \emph{weight}, given by a function $\omega: \Omega \to \N$. Given a coloring $\varphi: A \to \Omega$ of the objects, we define its \emph{total weight} to be the sum:
\begin{equation}
	|\varphi| = \sum_{x \in A}\omega \circ \varphi(x).
\label{eq:weight_sum}
\end{equation}
How many distinct colorings of $A$ with a given total weight are there, up to symmetries given by some group action $\Gamma_G$? Note that the total weight of any coloring in a given orbit is the same, as elements of $g$ merely permute the set $A$. Thus, the problem boils down to calculating the number of distinct orbits with a given total weight. P\'olya identified two necessary ingredients for a systematic answer to this question: generating functions, and an understanding of the cycle structure of elements of $G$ \cite{polya}.


\begin{dfnn}[Generating Function]
    Let $\omega: \Omega \to \N$ be an assignment of weights to some set $\Omega$. Suppose further that there are at most a finite number of elements of any given weight, that is, $|\omega^{-1}(n)|$ is finite for every $n \in \N$. Then the \emph{generating function} of $\omega$ is given by the polynomial:
    \begin{equation}
      \label{gen}
	    f_\omega(x) = \sum_{i = 0}^\infty |\omega^{-1}(i)|\ x^i.
    \end{equation}
\end{dfnn}

Generating functions are useful as they encode combinatorial data - in this case the number of colors of a given weight - as algebraic objects. In particular, we will need the following lemma:
\begin{lem}
\label{lemma:productcolor}
    Let $\omega_1: \Omega_1 \to \N$ and $\omega_2: \Omega_2 \to \N$ be assignments of weights to the sets $\Omega_1$ and $\Omega_2$ respectively. Define an assignment of weights to the set $\Omega_1 \times \Omega_2$ by $\omega: (x_1, x_2) \mapsto \omega_1(x_1) + \omega_2(x_2)$. Then $f_\omega(x) = f_{\omega_1}(x) \cdot f_{\omega_2}(x)$.
\end{lem}

Given a group action $\Gamma_G$ and an element $g \in G$, we denote by $C_k(g)$ the number of $k$-cycles in the unique disjoint cycle decomposition of $\Gamma_G(g)$. We can now encode information about the cycle structure of elements of $G$ in the following multivariate polynomial:
\begin{dfnn}[Cycle Index]
    Let $G$ be a finite group. Then the \emph{cycle index} of a group action $\Gamma_G$ on a finite set $A$ of cardinality $n$ is given by the polynomial \cite{brualdi}:
    \begin{equation}
	    Z_G(x_1, x_2, \ldots, x_n) = \frac1{|G|} \sum_{g \in G}{x_1^{C_1(g)}x_2^{C_2(g)}\cdots\ x_n^{C_n(g)}}.
\label{eq:pol}
    \end{equation}
\end{dfnn}

This cycle index will allow us to efficiently compute the number of distinct orbits of the group action. With this in mind, we are now in a position to state a version of the P\'olya counting theorem, answering the generalized problem given earlier:
\begin{teon}[P\'olya Enumeration Theorem]
Let $A$ be a finite set of objects, $\Omega$ a set of colors, $\omega: \Omega \to \N$ an assignment of weights to the colors with generating function $f_\omega$, and $\Gamma_G$ a group action of a finite group $G$ on $A$. Then $\Gamma_G$ passes naturally to a group action $\tilde \Gamma_G$ on $\Omega^{A}$, and a generating function by total weight for the number of distinct orbits of $\tilde \Gamma_G$ is given by:
\begin{equation}
	\text{Orbits}_{\tilde \Gamma_G}(x) = Z_G\left(f_w(x), f_w(x^2), \ldots, f_w(x^n)\right).
\label{eq:orb_num}
\end{equation}
\end{teon}

\subsection{P\'olya Enumeration Theorem for Bipartite Sets}
\label{sec:PETBS}

By considering multivariate generating functions, the P\'olya enumeration theorem can be generalized to the case where the colors take weights in $\N^k$. We will generalize the theorem in a different direction, however. Suppose we have a partition of $A$ into two parts, $A = X \sqcup Y$, and a group action $\Gamma_G$ on $A$. We would like to consider the problem of counting distinct colorings of $A$ under this symmetry, with the additional constraint that we color elements of $X$ from a set $\Omega_X$, and elements of $Y$ from a set $\Omega_Y$. To this end, we will say that a coloring $\varphi: A \to \Omega_X \sqcup \Omega_Y$ is \emph{valid} if $\varphi(x) \in \Omega_X \iff x \in X$ and $\varphi(x) \in \Omega_Y \iff x \in Y$.

There is an obstruction to this, however - the group action may map elements in $X$ to elements in $Y$ or vice versa. In this case, the extension of $\Gamma_G$ to the set of possible colorings is no longer well-defined, as there is no natural way to compare the sets of colors $\Omega_X$ and $\Omega_Y$. Fortunately, this is the only obstruction to proving a P\'olya-type theorem for this problem. This motivates the following definition:
\begin{dfnn} [Partition-Preserving Group Action]
Let $A = X \sqcup Y$, and let $\Gamma_G$ be a group action on $A$. Then we say that $\Gamma_G$ is \emph{partition-preserving} if for every $g \in G$, $gx \in X \iff x \in X$ and $gx \in Y \iff x \in Y$.
\end{dfnn}
The importance of this property is as follows. Suppose we have a group action $\Gamma_G$ on $A = X \sqcup Y$, and some element $g \in G$. Then $\Gamma_G(g)$ has a unique disjoint cycle decomposition given by $\Gamma_G(g) = C_1 \cdot C_2 \cdot \ldots \cdot C_k$. If $\Gamma_G$ is partition-preserving then each cycle $C_i$ is contained entirely in either $X$ or $Y$, and $\Gamma_G$ is in fact partition-preserving if and only if this is the case for every $g \in G$.

If $\Gamma_G$ is partition-preserving, then we define $C^X_k(g)$ to be the number of $k$-cycles in the disjoint cycle decomposition of $\Gamma_G(g)$ that are contained in $X$, and we define $C^Y_k(g)$ analogously. We will now define an analogue of the cycle index polynomial for the case of partition-preserving group actions. This will allow us to keep track of the cycle structure of elements of the group as well as which partition part each cycle acts on:
\begin{dfnn}[Bipartite Cycle Index]
	Let $G$ be a finite group and $A = X \sqcup Y$ a finite set of cardinality $n$. Then the \emph{bipartite cycle index} of a partition-preserving group action $\Gamma_G$ on $A$ is defined to be the polynomial:
	\begin{equation}
	    \tilde Z_G(x_1, \ldots, x_n, y_1, \ldots, y_n) = \frac 1{|G|} \sum_{g \in G} {x_1}^{C_1^X(g)} \cdots\ {x_n}^{C_n^X(g)} {y_1}^{C_1^Y(g)} \cdots\ {y_n}^{C_n^Y(g)}.
\label{eq:ZG}
    \end{equation}
\end{dfnn}

We can now generalize P\'olya's theorem to the case of partition-preserving group actions. We note that this theorem is used implicitly in \cite{polya} without proof.
\begin{teo}[Bipartite P\'olya Enumeration Theorem]
\label{theorem:bipartite_polya}
	Let $\Gamma_G$ be a partition preserving group action of a finite group $G$ on a finite set $A = X \sqcup Y$. Let $\Omega = \Omega_X \sqcup \Omega_Y$ be a set of colors, and let $\omega_X: \Omega_X \to \N^+$ and $\omega_Y: \Omega_Y \to \N^+$ be their assigned weights with respective generating functions $f_X$ and $f_Y$. If $\Phi$ is the set of valid colorings of $A$, then $\Gamma_G$ passes naturally to a group action $\tilde \Gamma_G$ on $\Phi$, and a generating function by total weight for the number of orbits of $\tilde \Gamma_G$ is given by:
	\begin{equation}
	    \text{Orbits}_{\tilde \Gamma_G}(x) = \tilde Z_G\left( f_X(x), \ldots, f_X(x^k), f_Y(x), \ldots, f_Y(x^k) \right).
\label{eq:orb_gx}
    \end{equation}
\end{teo}
\begin{proofn}
	We pass to a group action $\tilde \Gamma_G$ on $\Phi$ as follows. Given a valid coloring $\varphi \in \Phi$ and an element $g \in G$, we define the action of $g$ on $\varphi$ by $g\varphi: x \mapsto \varphi(gx)$. To compute a generating function for the number of orbits of $\tilde \Gamma_G$ by total weight, we will determine the generating functions for the number of fixed points of each $g \in G$ by total weight.

	Consider some $g \in G$. As $A$ is finite, there exists a unique disjoint cycle decomposition $\Gamma_G(g) = C_1 \cdot C_2 \cdot \ldots \cdot C_k$, where each $C_i$ is a cycle in the symmetric group $S_A$. Now suppose that $g$ fixes some valid coloring $\varphi \in \Phi$; that is, $g \varphi = \varphi$. Then, assuming  the cycle $C_i = (x_1\ x_2\ \ldots\ x_{k_i})$ for some $x_i \in A$, we have by definition that $\varphi(x_i) = (g\varphi)(x_i) = \varphi(gx_i) = \varphi(x_{i+1})$, and hence every element in the cycle must have the same color under $\varphi$. The number of colorings of $C_i$ that are fixed by $g$ is thus given by the generating function $f_X(x^{k_i})$ if $C_i$ lies in $X$, and $f_Y(x^{k_i})$ if $C_i$ lies in $Y$. We note that one of these two cases must occur for every cycle as $\Gamma_G$ is partition-preserving. By lemma \ref{lemma:productcolor}, then, the number of valid colorings of $A$ that are fixed by $g$ is given by the generating function:
    \begin{equation}
    	\text{Fix}_g(x) = f_X^{C_1^X(g)}(x) \cdots f_X^{C_k^X}(x^k) f_Y^{C_1^Y}(x) \cdots f_Y^{C_k^Y}(x^k).
	\end{equation}

	By Burnside's lemma, the number of orbits of $\tilde \Gamma_G$ of a particular weight is given by the average number of fixed colorings of that weight by elements $g \in G$. Applying Burnside's lemma for each possible weight, the number of orbits of $\tilde \Gamma_G$ is thus given by the generating function:
	\begin{align}
	    \text{Orbits}_{\tilde \Gamma_G}(x) &= \frac{1}{|G|}\sum_{g \in G} \text{Fix}_g(x) \nonumber \\
	    &= \frac{1}{|G|}\sum_{g \in G} f_X^{C_1^X(g)}(x) \cdots f_X^{C_k^X}(x^k) f_Y^{C_1^Y}(x) \cdots f_Y^{C_k^Y}(x^k) \nonumber \\
	    &= \tilde Z_G\left( f_X(x), \ldots, f_X(x^k), f_Y(x), \ldots, f_Y(x^k) \right).
    \end{align}
    \qed
\end{proofn}

We note that as a corollary of this proof, we can recover a bivariate generating function from this expression, where the coefficient of $x^ay^b$ represents the number of distinct colorings with total weight $a$ in $\Omega_X$, and total weight $b$ in $\Omega_Y$:
\begin{corn}
	A bivariate generating function by total weight in $\Omega_X$ and $\Omega_Y$, for the number of distinct colorings of $A$, is given by:
	\begin{equation}
    \label{eq:polya_cor}
	    \text{Orbits}_{\tilde \Gamma_G}(x, y) = \tilde Z_G\left( f_X(x), \ldots, f_X(x^k), f_Y(y), \ldots, f_Y(y^k) \right).
  \end{equation}
\end{corn}

\subsection{The Dihedral Group,  its Cycle Index and its Extension}
\label{sec:Dihedral}

To apply these results to the problem of counting distinct DNA necklaces, we will need to describe the relevant group action and compute its (bipartite) cycle index. The set of elements acted on by the group is given by the $\alpha$ containers in the DNA necklace and this set can be partitioned into two groups: containers of black beads and containers of white beads. We consider two DNA necklaces to be the same if one can be \emph{rotated} or \emph{reflected} into the other. These symmetries can be described by an action of the dihedral group, which we will denote by $D_{2M}$, where we have $\alpha = 2M$. {The rotational and reflective symmetries are what distinguishes the case of periodic DNA chains from linear, fixed boundary condition chains studied in \cite{Robin_Daudin_1999_JApplProb_36_179} and elsewhere.}

\begin{figure}[h!]
  \centering
  \includegraphics[scale=0.8]{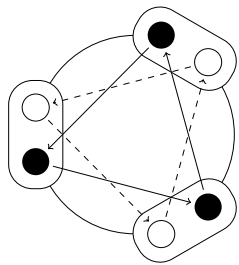}
  \caption{The action of a rotation $s \in D_{2M}$ on the DNA necklace.}
  \label{rot}
\end{figure}

\begin{figure}[h!]
  \centering
  \includegraphics[scale=0.8]{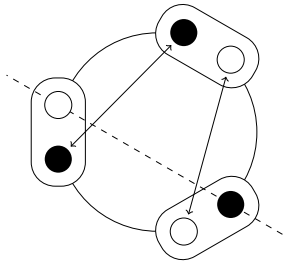}
  \caption{The action of a reflection $r \in D_{2M}$ on the DNA necklace, for the case where $M$ is odd.}
  \label{oddref}
\end{figure}

\begin{figure}[h!]
  \centering
  \includegraphics[scale=0.8]{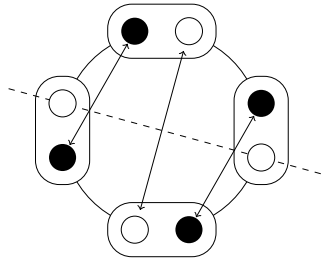}
  \caption{The action of a reflection $r \in D_{2M}$ on the DNA necklace, for the case where $M$ is even.}
  \label{evenref}
\end{figure}

A fundamental fact about $D_{2M}$ is that it is generated by two elements $r$ and $s$, where $r$ is a reflection satisfying $r^2 = 1$, and $s$ is a rotation of order $M$. Therefore, to describe a group action of $D_{2M}$ on a DNA necklace it suffices to give the action of $r$ and $s$. In Fig.~\ref{rot} the action of such a rotation on the necklace is illustrated, while in Figs.~\ref{oddref} and \ref{evenref} the action of a reflection is illustrated for the cases where $M$ is odd and even respectively. It is clear that the resulting group action is partition-preserving.

To compute the bipartite cycle index of this group action, we will treat reflections and rotations separately. To begin with, we can see from Fig. \ref{rot} that rotations act symmetrically on the black and white containers in the DNA necklace. Thus, the terms of the cycle index polynomial corresponding to rotations will be symmetric in the $x_i$ and $y_i$. The natural action of the cyclic group $C_M$ on the $M$ containers in a partition is given by \cite{vanlint}:
\begin{equation}
    Z_{C_M}(x_1, \dots, x_M) = \sum_{d \mid M} \varphi(d) x_d^{M/d},
\end{equation}
where $\varphi(d)$ is defined to be the number of natural numbers less that $d$ that are coprime to it (the \emph{Euler totient function}). Note that $1$ is considered to be coprime to all natural numbers, and so in particular $\varphi(d) > 0$. Exactly half of the elements of $D_{2M}$ are rotations, and thus the rotational part of the bipartite cycle index $\tilde Z_{D_{2M}}$ is given by $\frac 12 \sum_{d \mid M} \varphi(d) x_d^{M/d} y_d^{M/d}$.

The reflective part of the group $D_{2M}$, on the other hand, acts differently depending on the parity of $M$. Suppose first that $M$ is odd, in which case a typical reflection is illustrated in Fig. \ref{oddref}. Each of the $M$ possible reflections occur across an axis consisting of one black container and one white container, both of which are fixed by the reflection. The rest of the containers are split into $2$-cycles, and thus the bipartite cycle index $\tilde Z_{D_{2M}}$ for odd $M$ is given by:
\begin{equation}
  \label{eq:M-odd}
	\tilde Z_{D_{2M}}(x_1, \ldots, x_M, y_1, \ldots, y_M) = \frac 12 \sum_{d \mid M} \varphi(d) x_d^{M/d} y_d^{M/d} + \frac 12 x_1y_1x_2^{(M-1)/2}y_2^{(M-1)/2}.
\end{equation}
If $M$ is even, a typical reflection is illustrated in Fig. \ref{evenref}. In this case, each possible reflection occurs across an axis consisting of either two white containers or two black containers. The rest of the containers again split into $2$-cycles. Thus the bipartite cycle index $\tilde Z_{D_{2M}}$ for even $M$ is given by:
\begin{equation}
  \label{eq:M-even}
	\tilde Z_{D_{2M}}(x_1, \ldots, x_M, y_1, \ldots, y_M) = \frac 12 \sum_{d \mid M} \varphi(d) x_d^{M/d} y_d^{M/d} + \frac 14 x_1^2x_2^{(M-2)/2}y_2^{M/2} + \frac 14 y_1^2y_2^{(M-2)/2}x_2^{M/2}.
\end{equation}

\subsection{Generating Functions as Formal Power Series}
\label{sec:GenFun}

In our particular application of P\'olya theory, the \emph{elements} we are coloring are the $\alpha$ containers in the DNA necklace and the \emph{color} of a particular container is defined to be the number of black or white beads it contains. As each container must contain at least one bead, the set of colors is given by $\N^+$. We are interested in the total number of black and white beads, so the weight of each color will be given quite simply by $\omega(n) = n$ for each $n \in \N^+$. This weighting corresponds to the generating function (\ref{gen}) $f_\omega(x) = x + x^2 + x^3 + \cdots$.

To compute the number of distinct DNA necklaces with $N_{AT}$ white beads and $N_{GC}$ black beads, we need to calculate the coefficient of $x^{N_{AT}}y^{N_{GC}}$ in (\ref{eq:polya_cor}), where the bivariate cycle index is given by the appropriate $\tilde Z(D_{2M})$ from Sect.~\ref{sec:Dihedral} and the weight generating function is given by $f_\omega(x)$. This requires us to calculate the coefficients of specific terms in $f_\omega(x)^n = (x + x^2 + x^3 + \dots)^n$ for potentially large $n$. However, doing this expansion naively requires many computing steps, whose number grows exponentially fast as $n$ increases. Thus, this approach is impractical.  Fortunately, there exists a way to bypass this problem: treating $f_\omega(x)$ as a formal power series, we can manipulate it into a form that makes such computations significantly faster.

An introduction to the theory of formal power series can be found, for instance, in \cite{zariski}. For our purposes, we will only need the fact that a form of the binomial theorem holds in this setting:
\begin{lem}
\label{lemma:binomial}
   Letting $(1 - x)^{-n}$ denote the formal inverse of $(1 - x)^n$, we have:
  \begin{equation}
    (1-x)^{-n} = \sum_{k=0}^{\infty} \binom{n+k-1}{n-1} x^k.
    \label{eq:power}
  \end{equation}
\end{lem}
This implies the following useful lemma regarding powers of $f_\omega(x)$:
\begin{lem}
	As a formal power series $f_\omega(x)^n$ can be written as $f_\omega(x)^n = \sum_{k=0}^{\infty} \binom{n+k-1}{n-1}x^{n + k}$.
\end{lem}
\begin{proofn}
	Note that $x f_\omega(x) = x^2 + x^3 + \dots = f_\omega(x) - x$. Rearranging this for $f_\omega(x)$, we see that $f_\omega(x) = x(1 - x)^{-1}$, and hence $f_\omega(x)^n = x^n(1 - x)^{-n}$. The result now follows from lemma \ref{lemma:binomial}.
\end{proofn}
In contrast to naively expanding powers of $f_\omega(x)$, computing binomial coefficients is computationally inexpensive, taking at most a linear number of steps in $n$. 

We now list a few results that will come in handy later, when we describe an explicit algorithm for computing the number of distinct DNA necklaces with the given constraints.
\begin{lem}
	The coefficient of $x^r$ in $f_\omega(x^a)^b$ is given by:
	\begin{equation}
    \label{eq:coeff_mon}
	    \left[ f_\omega(x^a)^b \right]_r =
	        \begin{cases}
	            1                      &\quad\text{if }b=0 \text{ and } a=0 \\
	            0                      &\quad\text{if }b=0 \text{ and } a>0 \\
	            0                      &\quad\text{if }b>0 \text{ and }a \nmid r \text{ or } r < ab \\
	            \binom{r/a - 1}{b - 1} &\quad\text{otherwise.}
                         \end{cases}
    \end{equation}
\end{lem}
\begin{lem}
	The coefficient of $x^r$ in $f_\omega(x^{a_1})^{b_1} \cdot f_\omega(x^{a_2})^{b_2}$ is given by:
	\begin{equation}
    \label{eq:coeff_bin}
	    \left[ f_\omega(x^{a_1})^{b_1} \cdot f_\omega(x^{a_2})^{b_2} \right]_r = \sum_{k=0}^r \left[ f_\omega(x^{a_1})^{b_1} \right]_k \left[ f_\omega(x^{a_2})^{b_2} \right]_{r - k}.
    \end{equation}
\end{lem}

\section{The Algorithm for Computing the Number of Distinct Valid Necklaces}
\label{sec:algorithm}

Now we are able to evaluate the number of distinct necklaces, which correspond to a particular value of alternations $\alpha$. The algorithm is fairly straightforward and efficient. Its implementation requires the following steps:
\begin{enumerate}[a)]
  \item Set constraint parameters, $N_{AT}$, $N_{GC}$, and $\alpha = 2M$.
  \item Choose partitioned cycle index polynomial of the Dihedral group based on parity of $M$. If $M$ is odd, use (\ref{eq:M-odd}), while for  $M$ even use (\ref{eq:M-even}).
  \item By the corollary to P\'olya's Enumeration Theorem (\ref{eq:polya_cor}), we know that the number of necklaces, up to symmetry, is given by
  \begin{align}
    \label{eq:polya3}
    \text{Orbits}_{\tilde \Gamma_G}(x, y) =\  &\tilde Z_G\left( f_X(x), \dots, f_X(x^k), f_Y(y), \dots, f_Y(y^k) \right).
  \intertext{If $M$ is odd using the outcome of the previous step we get}
    \label{eq:expl_zno}
    \text{Orbits}_{\tilde \Gamma_G}(x, y) = \ &\frac{1}{2M}\sum_{d|M} \varphi(d)f^{M/d}(x^d)f^{M/d}(y^d) \nonumber\\
    &+ \frac{1}{2}f(x)f(y)f^{(M-1)/2}(x^2)f^{(M-1)/2}(y^2).
  \intertext{If $M$ is even, then we have}
    \label{eq:expl_zne}
    \text{Orbits}_{\tilde \Gamma_G}(x, y) = \ &\frac{1}{2M}\sum_{d|M} \varphi(d)f^{M/d}(x^d)f^{M/d}(y^d)\nonumber\\
     &+ \frac{1}{4}f^2(x)f^{(M-2)/2}(x^2)f^{M/2}(y^2) + \frac{1}{4} f^2(y)f^{(M-2)/2}(y^2)f^{M/2}(x^2).
  \end{align}
  \item Every term in the polynomial produced by (\ref{eq:polya3}) will be of the form in (\ref{eq:coeff_mon}) or (\ref{eq:coeff_bin}). The number of necklaces with $N_{AT}$ white beads and $N_{GC}$ black beads is given by the coefficient of the term $x^{N_{AT}}y^{N_{GC}}$. To calculate the total number of necklaces, simply sum over each of these terms appearing in the polynomial.
\end{enumerate}
A Python computer code  implementating this algorithm is presented in the Appendix.

In order to illustrate the application of this algorithm let us consider a simple, but not trivial case: We set  $\alpha = 2M = 10$, $N_\text{AT}=8$, $N_\text{GC}=6$. Clearly $M=5$ is odd, so identifying white beads with AT base pairs and black beads with GC base pairs, we have the cycle index
\begin{align}
  \tilde{Z}(\tilde{D}_{10}) &= \frac{1}{2}Z(\tilde{C}_5) + \frac{1}{2} x_1 y_1(x_2)^{2}(y_2)^{2} \nonumber\\
  & = {\frac  {1}{5}}\sum _{{d|5}}\varphi (d)(x_{d})^{{5/d}}(y_{d})^{{5/d}} + \frac{1}{2} x_1 y_1(x_2)^{2}(y_2)^{2}.
\end{align}
Now the partitioned P\'olya Enumeration Theorem tells us that we can put the generating functions $f_W\left(x^d\right)$ and $f_B\left(y^d\right)$ in place of the $x_d$ and $y_d$ respectively to find the generating function of fixed orbits. So we have
\begin{align}
   \text{Orbits}_{\tilde \Gamma_G}(x, y) = &\frac{1}{2\cdot5}\left[ 1 (x+x^2 + x^3+\dots)^5(y+y^2+y^3+\dots)^5\nonumber \right.\\
   + &\left.4(x^5 + x^{10} + x^{15}+\dots)(y^5+y^{10}+y^{15}+\dots) \right] \nonumber\\
  +&\frac{1}{2}(x + x^2 + \dots)(x^2+x^4 + \dots )^2(y + y^2 + \dots)(y^2 + y^4 + \dots)^2.
\end{align}
Let us first look at the cyclic part. Since 5 is prime, the only two integers that divide it are 1 and 5, so this polynomial will be
\begin{equation*}
  \frac{1}{2\cdot5}\left[ 1 (x+x^2 + x^3+\ldots)^5(y+y^2+y^3+\ldots)^5 + 4(x^5 + x^{10} + x^{15}+\ldots)(y^5+y^{10}+y^{15}+\ldots)\right].
\end{equation*}
Now we try to extract the coefficients of terms that are allowed. These are the terms in $x^{N_\text{AT}}$ and $y^{N_\text{GC}}$ and we can use (\ref{eq:coeff_mon}) in order to  calculate these coefficients directly. In this case, there will be no contribution from the second term, as there are no terms in $x^8$ and $y^6$. So the total cyclic contribution will be (with $r=8$ and $r=6$ for the respective cases and $a=1, \ b=5$ for both)
\begin{equation*}
  \frac{1}{10} \binom{N_{\text{GC}}-1}{5-1} \binom{N_{\text{AT}}-1}{5-1}= \frac{1}{10}\binom{5}{4}\binom{7}{4} = \frac{175}{10}.
\end{equation*}
Then the same coefficient identifying process can be followed for the reflective part. Now the polynomial is given by
\begin{equation*}
  \frac{1}{2}(x + x^2 + \ldots)(x^2+x^4 + \ldots )^2(y + y^2 + \ldots)(y^2 + y^4 + \ldots)^2.
\end{equation*}
So for both $x$ and $y$ the coefficients will come from the product of two series, one of them squared. Thus,  the relevant terms will come in a series of products given in (\ref{eq:coeff_bin}). In $y$ the sum of coefficients contracts to a single element. That contribution is simply $\binom{1}{0}\binom{1}{1}=1$. In $x$ however, there will be terms from $x^2\cdot x^6$ as well as $x^4\cdot x^4$. So then, the sum will be
\begin{equation*}
  \binom{1}{0}\binom{3}{1} + \binom{3}{0}\binom{1}{1} = 4,
\end{equation*}
giving a total contribution of $\frac{1}{2}(1+4) + \frac{175}{10} = 20$. Thus there are 20 DNA chains with 8 AT base pairs, 6 GC base pairs and 10 alternations.

\section{Numerical Results}
\label{sec:NumRes}

\begin{figure}[h!]
  \includegraphics[width=0.47\textwidth]{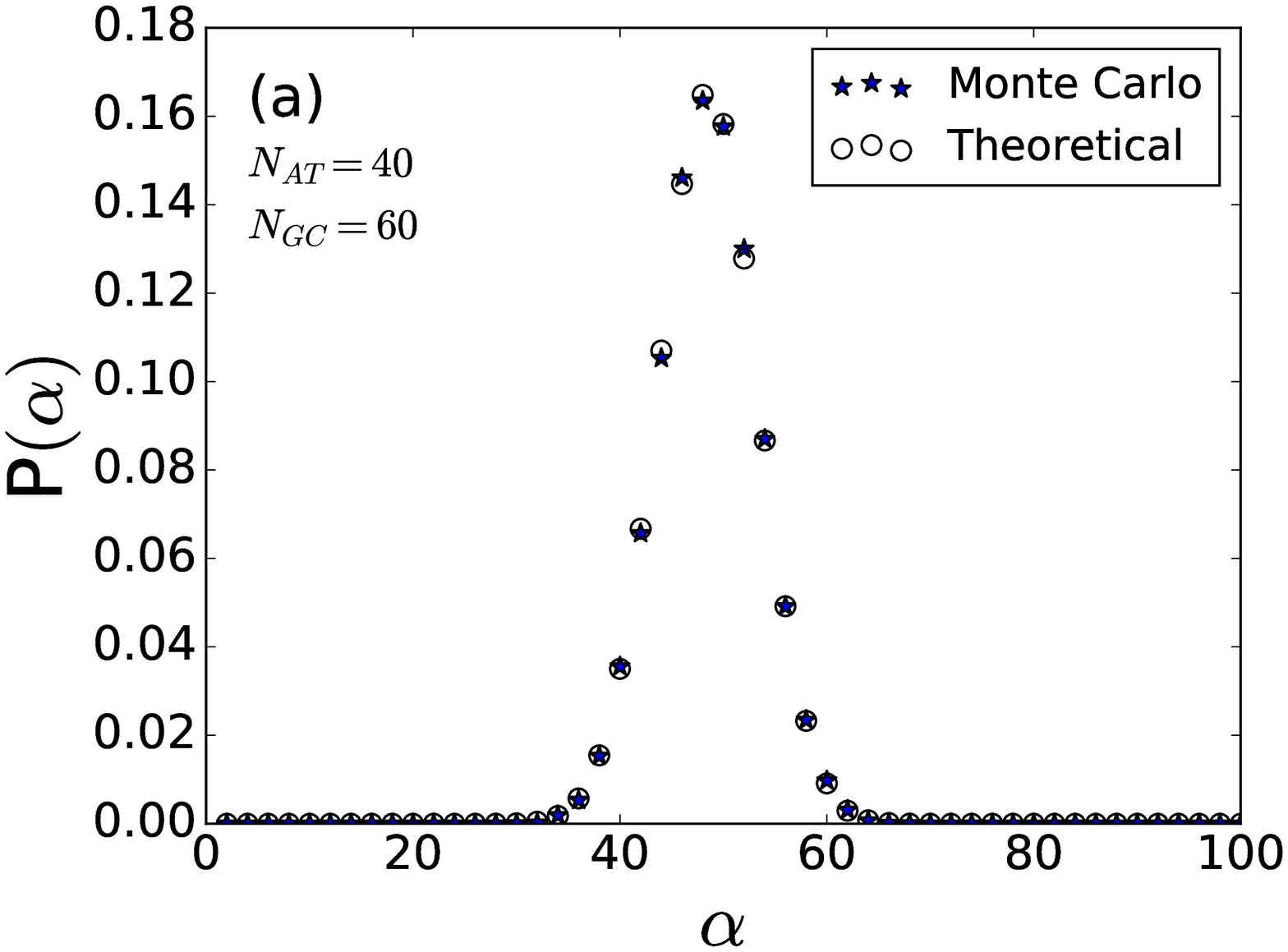}
  \includegraphics[width=0.47\textwidth]{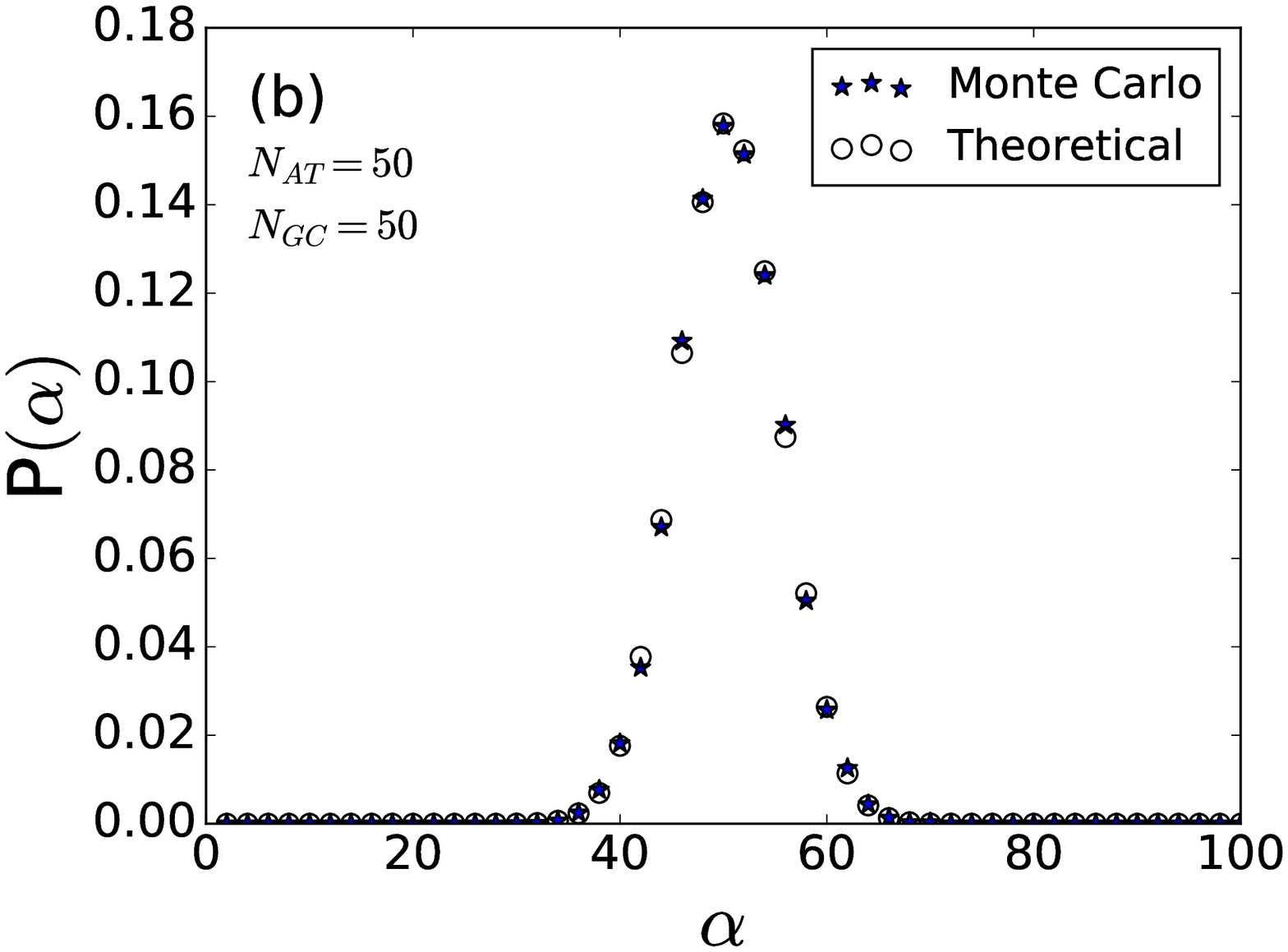}
  \includegraphics[width=0.5\textwidth]{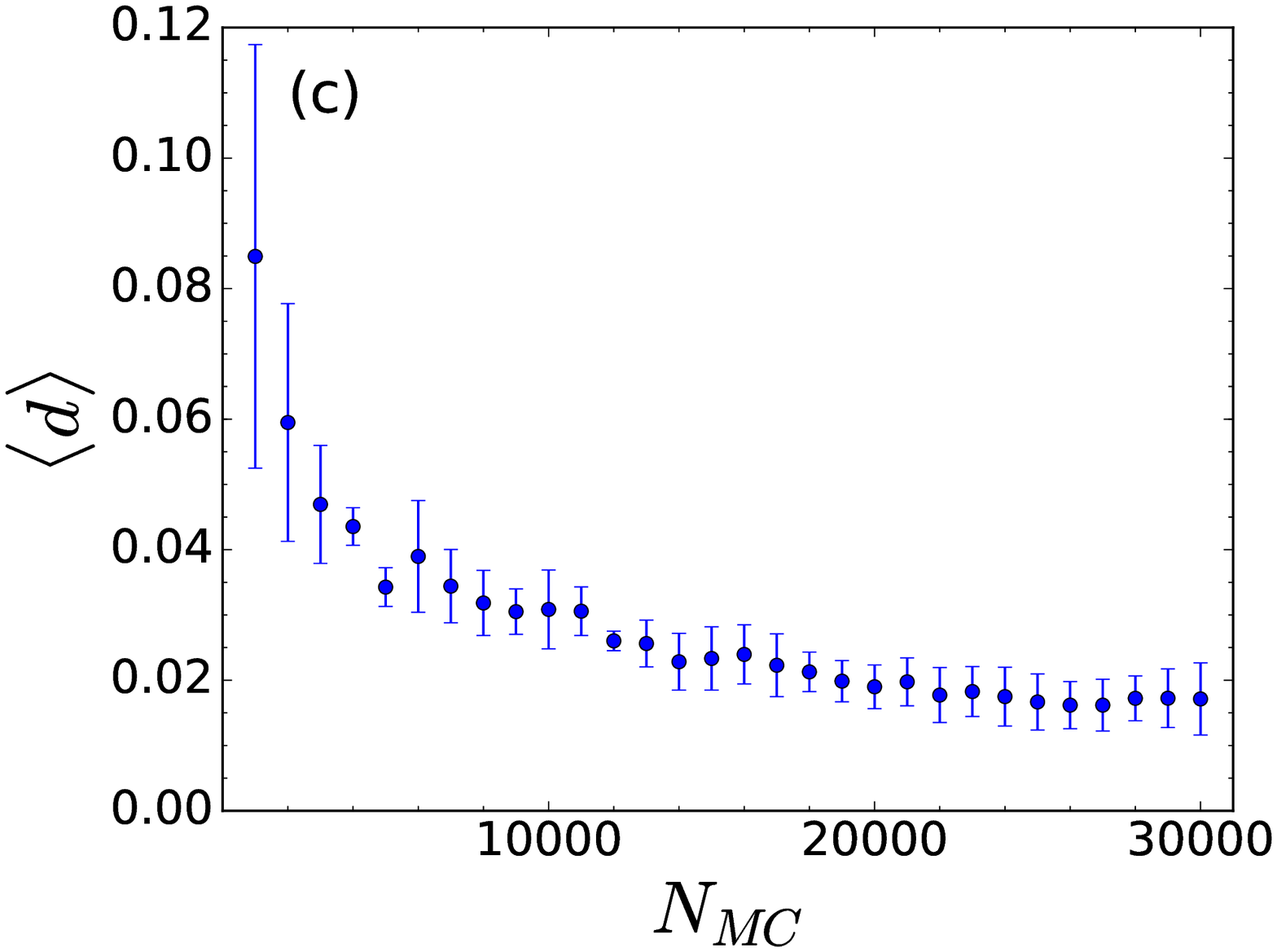}
  \caption{Comparison of the pdf $P(\alpha)$ of the number of alternations $\alpha$, obtained by the algorithm presented in Sect.~\ref{sec:algorithm} [empty circles in panels (a) and (b)] and by randomly created DNA chains of $N=100$ base pairs through MC simulations  [filled stars in panels (a) and (b)]. The pdfs for $N_{AT}=40$, $N_{GC}=60$ and $N_{AT}=50$, $N_{GC}=50$  are presented in panels (a) and (b) respectively. The number of MC simulations used in (a) and (b) are $N_{MC}=20000$.
  (c) The evolution of the average total absolute difference $\langle d \rangle$ between the theoretically and the numerically obtained pdfs as a function of $N_{MC}$ for the case of $N_{AT}=50$, $N_{GC}=50$. The values of $\langle d \rangle$ are obtained as the average of the quantity (\ref{eq:d}) evaluated for 5 different sets of $N_{MC}$ runs. The error bars denote the corresponding standard deviations.}
  \label{ds}
\end{figure}

The developed  algorithm for calculating the number  of distinct DNA chains having $\alpha$ alternations can be used to produce the pdf of $\alpha$, $P(\alpha)$,  which afterwards can be compared to pdfs numerically obtained from Monte-Carlo (MC) simulations. In Figs.~\ref{ds}(a) and (b) we present such pdfs for a DNA chain containing $N=100$ base pairs. In particular, we consider the case of $N_{AT}=40$, $N_{GC}=60$ in Fig.~\ref{ds}(a) and the case of  $N_{AT}=50$, $N_{GC}=50$ in Fig.~\ref{ds}(b).
From Figs.~\ref{ds}(a) and (b) we clearly see that the results obtained by the  algorithm presented in Sect.~\ref{sec:algorithm} (empty circles) and by MC simulations of DNA chains with $N=100$ base pairs (filled stars) agree very well. The slight differences between them are to be expected, as the number of possible chains is generally very large. For instance, in the case of $N_{AT} = 50$, $N_{GC}=50$ and $\alpha=50$, the number of possible DNA chains is of the order of $10^{25}$ possible necklaces. Thus, in general, the  number of performed MC simulations cannot  get close to the actual total number of possible chains. Nevertheless, although the results of Figs.~\ref{ds}(a) and (b) were obtained by  only $N_{MC}=20000$ MC simulations they manage to capture the theoretically obtained pdf quite accurately. Of course it is expected that increasing the number of MC simulations will improve the accuracy of the numerical results.
As a measure of this accuracy we can consider the total absolute difference
\begin{equation}
  d(N_{MC}) = \sum_{\alpha} |P_{MC}(N_{MC}, \alpha)-P(\alpha)|,
  \label{eq:d}
\end{equation}
between the two distributions. In (\ref{eq:d}) $P_{MC}(N_{MC}, \alpha)$ is the probability of $\alpha$ alternations obtained by $N_{MC}$ MC simulations, $P(\alpha)$ is the one obtained theoretically, while the sum is performed over all possible values of $\alpha$. From the results of Fig.~\ref{ds}(c) where we plot the averaged value of $d(N_{MC})$ over 5 sets of $N_{MC}$ MC simulations as a function of $N_{MC}$ we  see that as the number of simulations increases, the numerical results get closer to the theoretical ones.

The results of Fig.~\ref{ds} clearly show that in order to study the dynamical properties of DNA chains, statistical analysis performed over a few thousands of MC generated random chains (even of the order of 5000) would suffice, as such numbers of MC simulations are enough for capturing quite accurately the influence of alternations on the system's dynamics.

\begin{figure}[h!]
  \centering
  \includegraphics[width=0.57\textwidth]{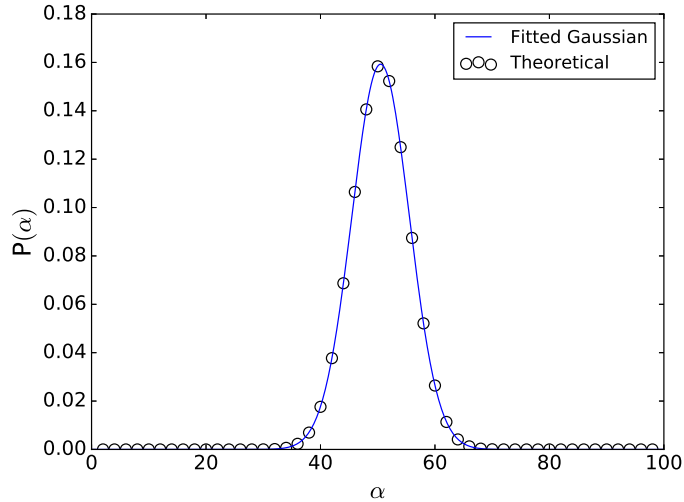}
  \caption{Fitting by a Gaussian of the theoretical pdf of Fig.~\ref{ds}(b) (empty circles) with $N_{AT} = 50$, $N_{GC}=50$. The mean of the  Gaussian is $\alpha_0=50.5$ and standard deviation $\sigma_{\alpha}=5.1$.}
  \label{fit}
\end{figure}

\begin{figure}[h!]
  \centering
  \includegraphics[width=0.6\textwidth]{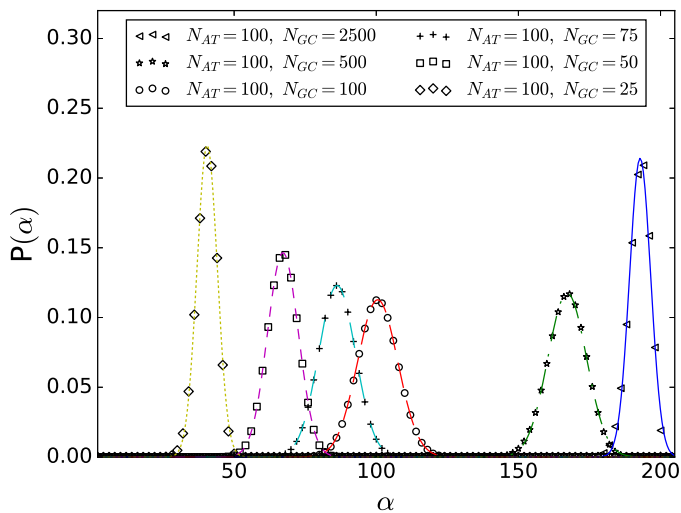}
  \caption{Pdfs of $\alpha$ for fixed number of AT base pairs ($N_{AT}=100$) and increasing values of $N_{GC}$. Points correspond to the theoretically obtained values of the pdfs, while curves correspond to the Gaussian fits of these points. Note that even for long DNA chains the value of $\alpha$ cannot exceed  $\alpha=200$.}
  \label{lop}
\end{figure}

The shape of the pdfs in Figs.~\ref{ds}(a) and (b) suggests that they could possibly be fitted by Gaussian distributions. This is actually true as we can see from the results of Fig.~\ref{fit}, where we performed such a fit for the theoretically obtained pdf of Fig.~\ref{ds}(b).
The Gaussian approximation of the pdfs has several advantages as it allows us to easily quantify the influence of different variables on the number of alternations.
Let us first look at the effect of increasing the number of only one type of base pair, keeping constant the number of the other type of base pair. In Fig.~\ref{lop} we present some pdfs of $\alpha$ for $N_{AT}=100$ and increasing values of $N_{GC}$
from 25 up to 2500.
Starting from small values of $N_{GC}$, we find a very ``lopsided'' and narrow distribution which as $N_{GC}$ increases becomes gradually more symmetric and spreads out, up to a value of $N_{GC}=200$. Then, increasing $N_{GC}$ further, as the numbers of different types of base pairs become more dissimilar we again find gradually more unbalanced pdfs with sharp peaks. The very ``lopsided'' base pair distributions are obtained when the minority base pairs are significantly less than the majority ones and therefore are spread out and isolated among the others. In this case the distribution is sharply peaked around the corresponding maximum possible number of alternations. For the $N_{AT}=100$, $N_{GC}=25$ case this number is $\alpha=50$, while for the $N_{AT}=100$, $N_{GC}=2500$ case it is $\alpha=200$.
\begin{figure}[h!]
  \centering
  \includegraphics[scale=0.26]{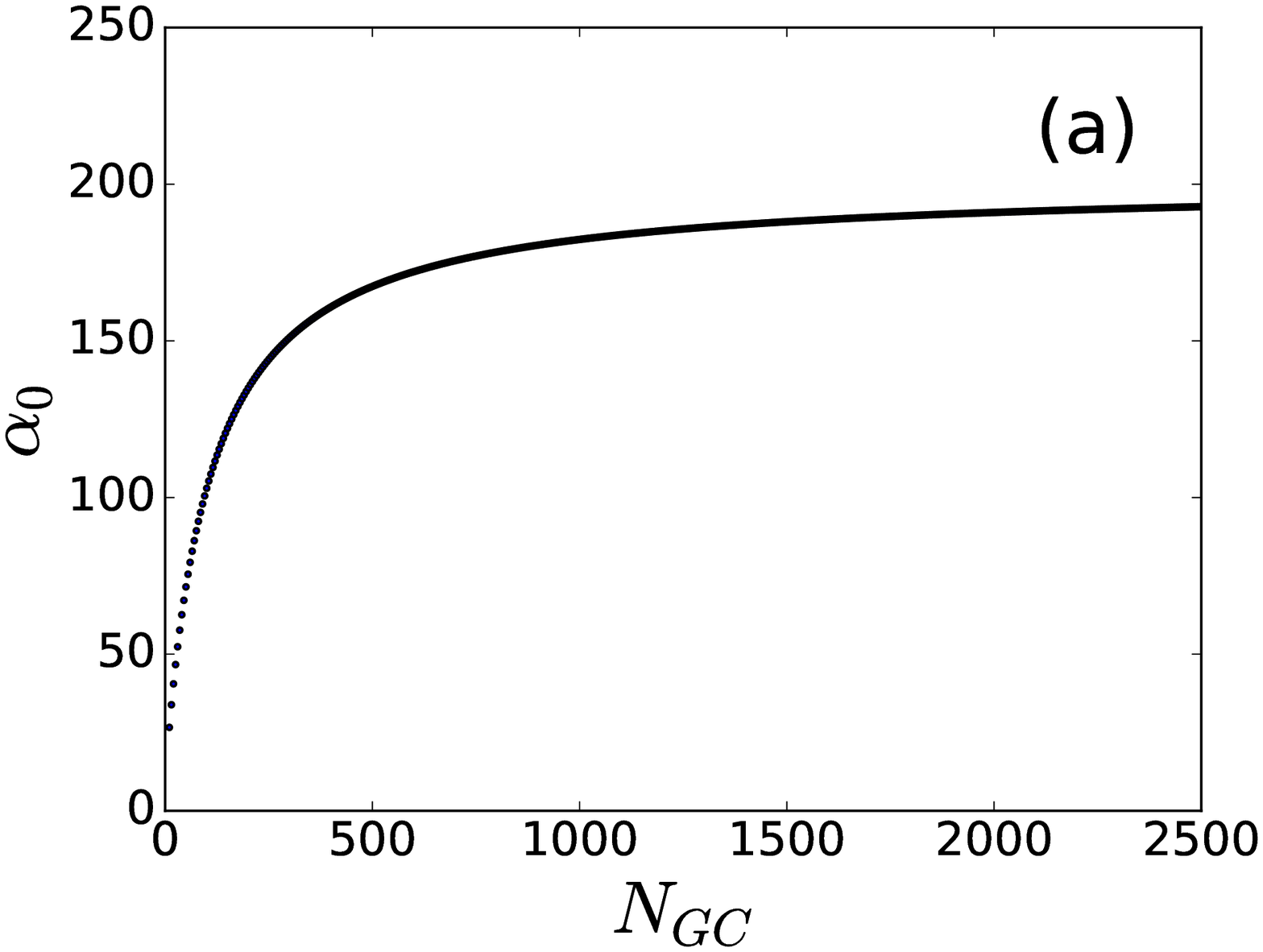}
  \includegraphics[scale=0.26]{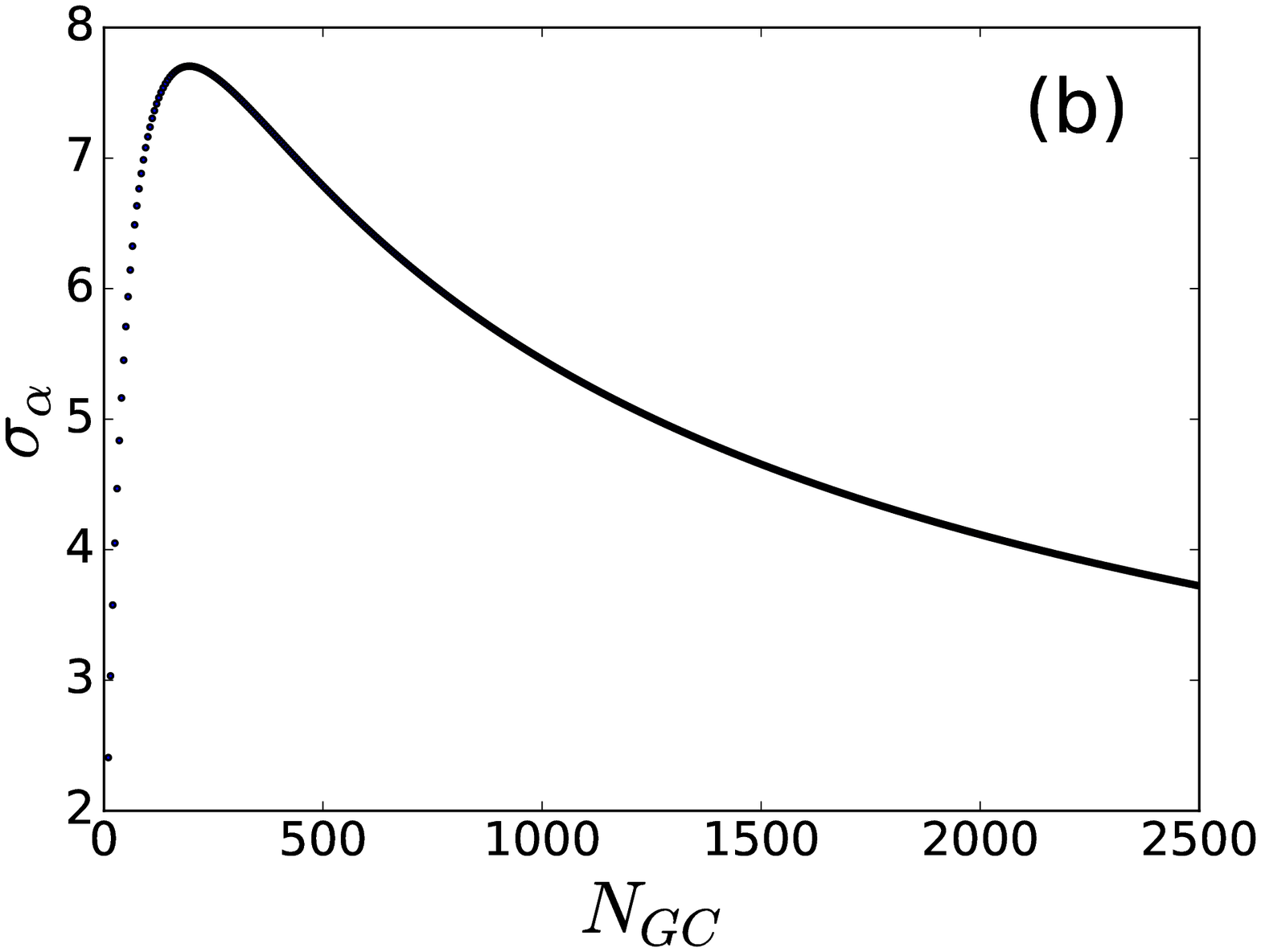}
  \includegraphics[scale=0.26]{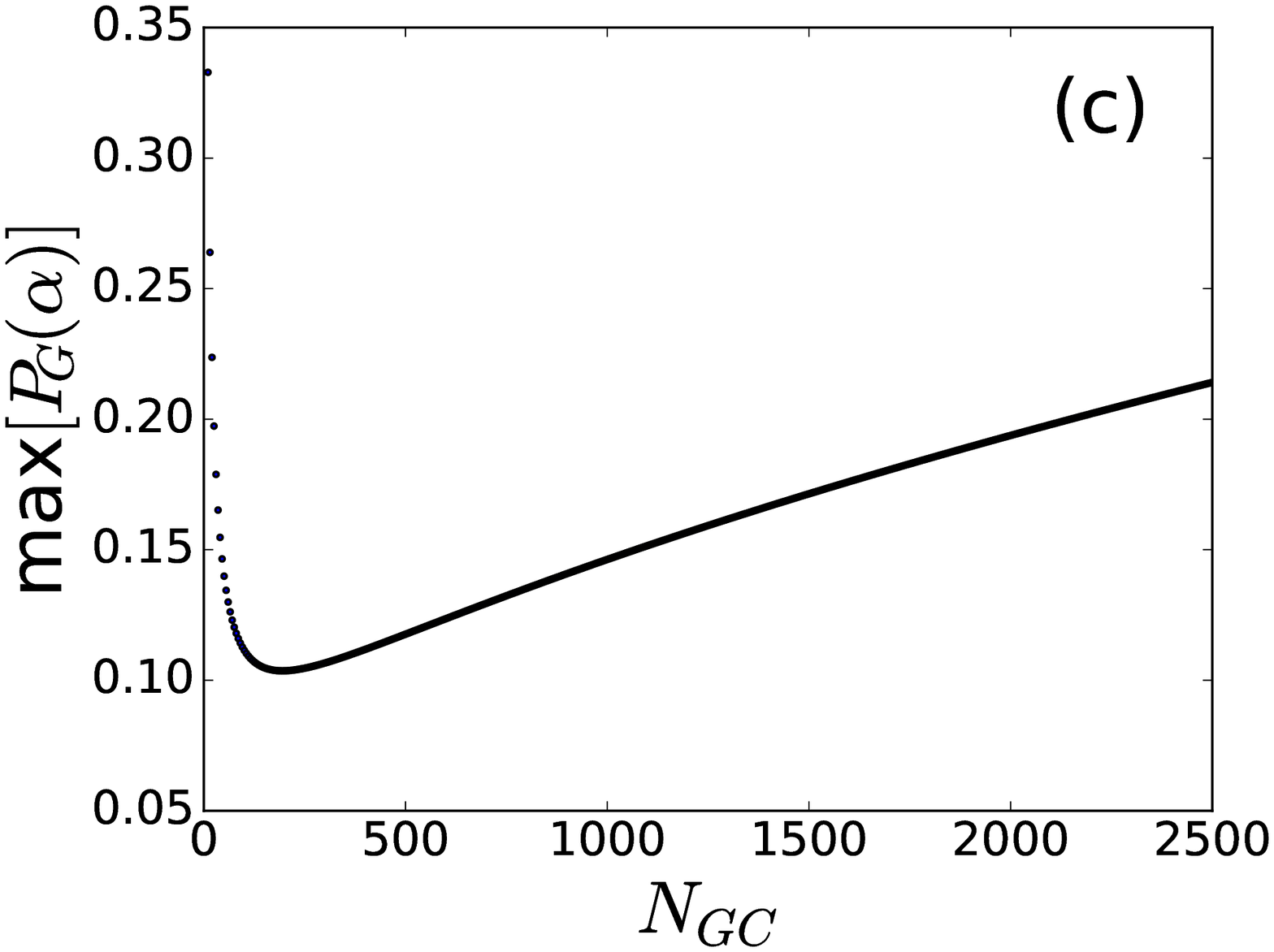}
  \caption{The effect of increasing the number $N_{GC}$ of the GC base pairs for a fixed number of AT base pairs ($N_{AT}=100$) on the Gaussian fit $P_G(\alpha)$ of the pdf values of $\alpha$, and in particular  on (a) the mean value $\alpha_0$, (b) the standard deviation $\sigma_{\alpha}$ and (c) the maximum probability $\max \left[  P_G(\alpha) \right]$. Some of these pdfs are shown in Fig.~\ref{lop}.}
  \label{maxes}
\end{figure}

These changes of the distributions are quantitatively presented in Fig.~\ref{maxes} through the variations of the fitted Gaussian characteristics. The increase of the mean value $\alpha_0$ of the Gaussian fits as the number $N_{GC}$ increases is shown in Fig.~\ref{maxes}(a). The upper limit of $\alpha_0$ is 200, when $N_{GC}$ becomes much larger than $N_{AT}$.
The dependence of the width (standard deviation) $\sigma_{\alpha}$ of the Gaussian fits on $N_{GC}$ is depicted in Fig.~\ref{maxes}(b). The initial increase with $N_{GC}$ corresponds to the spreading out of the distributions when the numbers of base pairs become more similar. Further increase of the $N_{GC}$ values pushes the pdfs to the other extreme and the lopsidedness comes through again, resulting in narrower distributions (see Fig.~\ref{lop}). This results in the decrease of $\sigma_{\alpha}$ for large values of $N_{GC}$. Finally in Fig.~\ref{maxes}(c) we observe that as $N_{GC}$ increases the maximum probability of the pdfs initially decreases rapidly and then increases slowly, in accordance with the results of Fig.~\ref{lop} and of course with the fact that it is inversely proportional to the standard deviation of the Gaussian fit.

\begin{figure}[h!]
  \centering
  \includegraphics[width=0.355\textwidth]{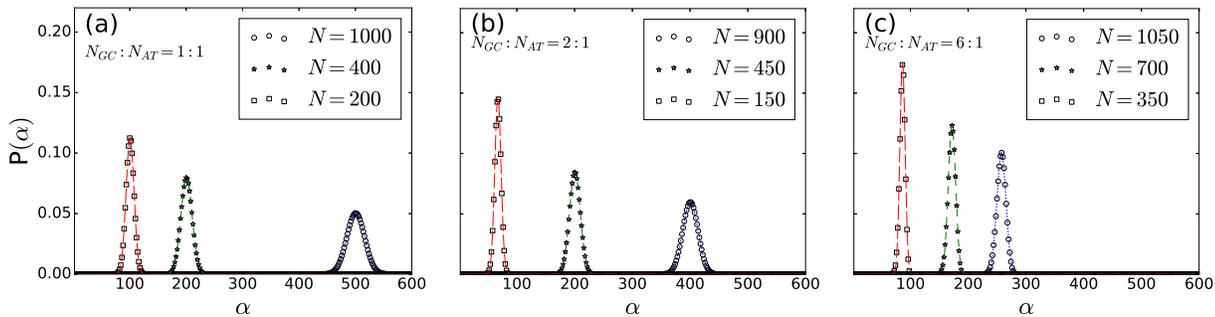}
  \includegraphics[width=0.31\textwidth]{Alpha_2to1no_lab76b.eps}
  \includegraphics[width=0.31\textwidth]{Alpha_6to1no_lab76b.eps}
  \caption{Pdfs of $\alpha$ for fixed ratios $N_{GC}:N_{AT}=$ $1:1$ (a), $2:1$ (b) and $6:1$ (c). Points correspond to the theoretically obtained values of the pdfs, while curves correspond to the Gaussian fits of these points.}
  \label{bal}
\end{figure}

\begin{figure}[h!]
  \includegraphics[width=0.31\textwidth]{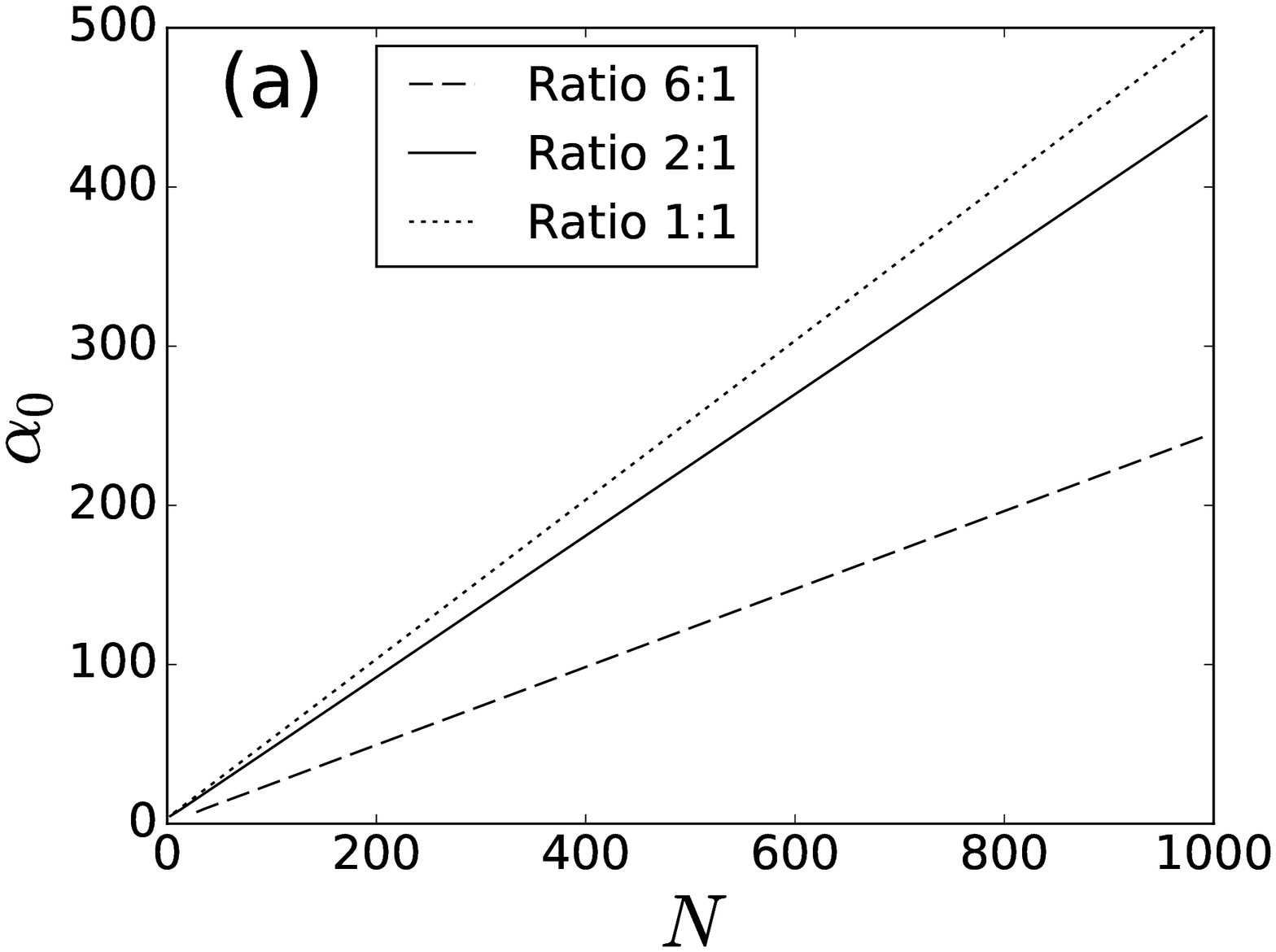}
  \includegraphics[width=0.31\textwidth]{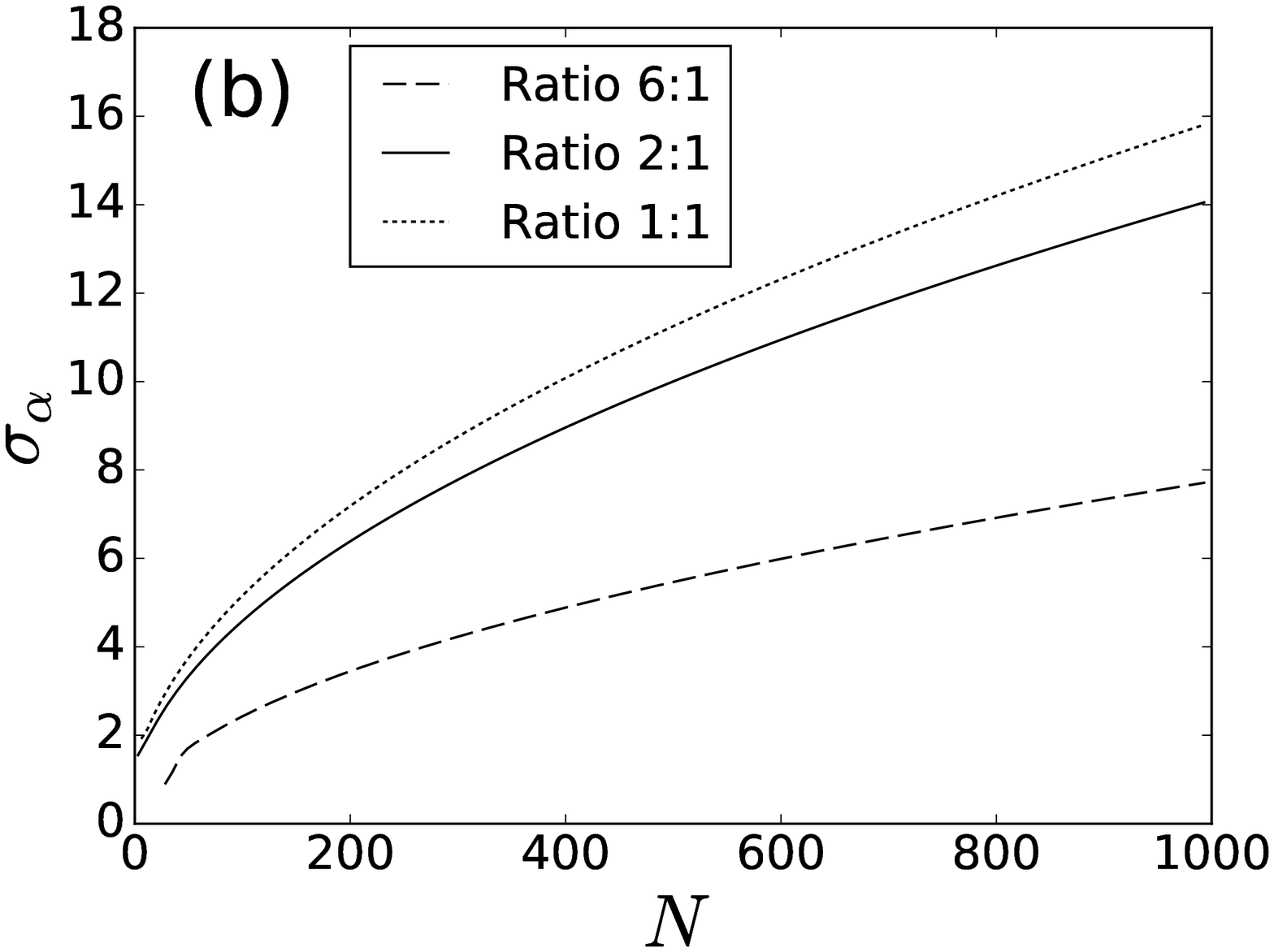}
  \includegraphics[width=0.31\textwidth]{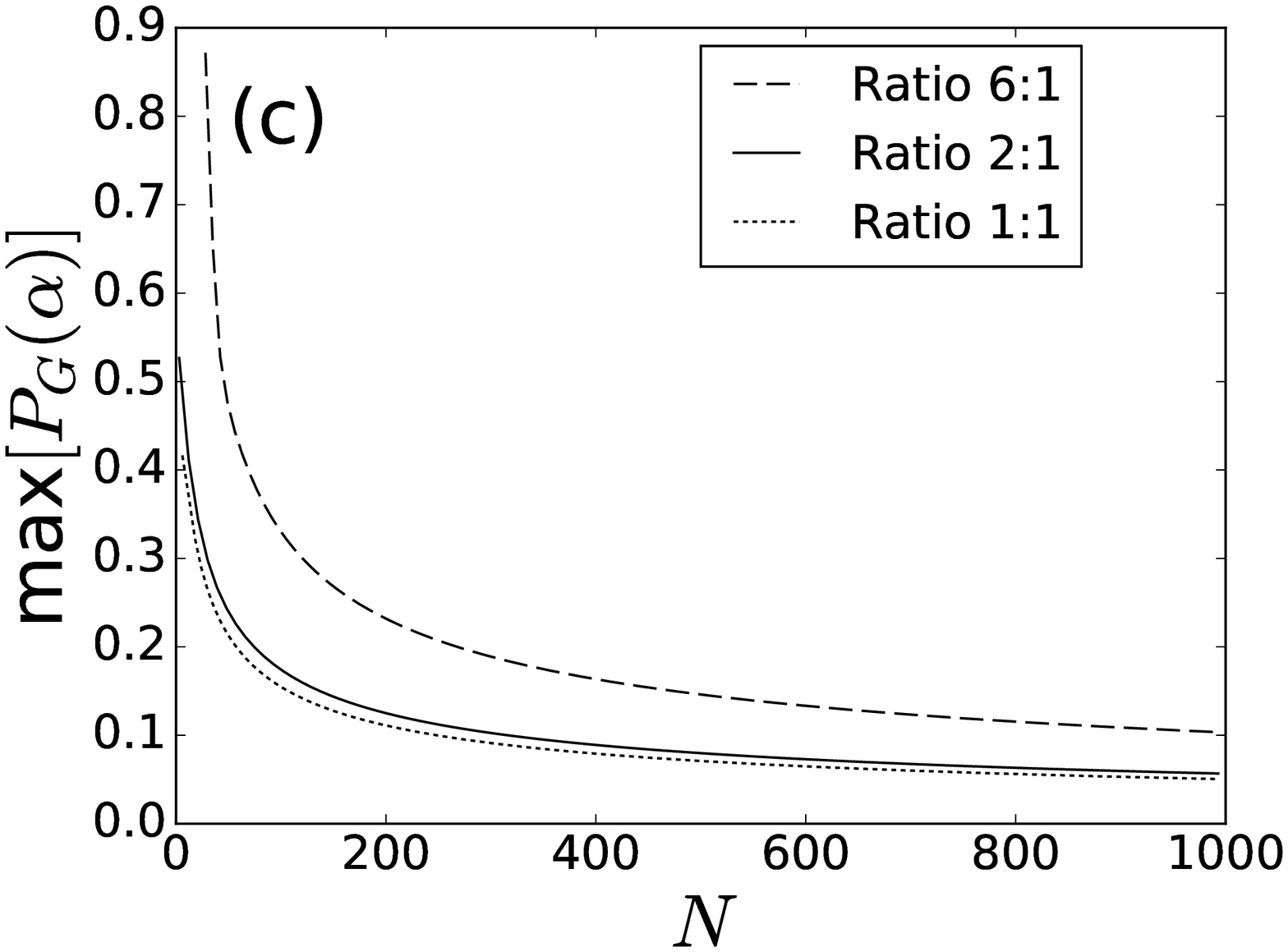}
  \caption{The effect of increasing the total number of base pairs $N$ for fixed ratios $N_{GC}:N_{AT}$ on the parameters of the Gaussian fit $P_G(\alpha)$ of the pdf for $\alpha$: (a) the mean value $\alpha_0$, (b) the standard deviation $\sigma_{\alpha}$ and (c) the maximum probability $\max \left[  P_G(\alpha) \right]$. Some of these pdfs are shown in Fig.~\ref{bal}.}
  \label{rats}
\end{figure}

Let us now focus our attention on the effect of the increment of the total number of base pairs $N =N_{AT}+N_{GC}$, i.e.~the total `length' of the DNA chain, when the ratio $N_{GC}:N_{AT}$ is kept constant. Such cases are presented in Fig.~\ref{bal}, where we plot several pdfs for different values of $N$ but for fixed ratios $N_{GC}:N_{AT}$.
In particular,  the values of the ratios $N_{GC}:N_{AT}$ are  $1:1$ in panel (a), $2:1$ in (b) and $6:1$ in (c). In all cases the pdfs are fitted by appropriate Gaussian distributions whose characteristics are plotted in Fig.~\ref{rats} as a function of $N$.
From the results of Figs.~\ref{bal} and \ref{rats} we see that as the total number $N$ of base pairs increases the pdfs become more broad, and consequently their maximum value decreases. This means that for large $N$ more $\alpha$ values have a relatively high probability to appear in a randomly created DNA chain. In addition, increasing the ratio $N_{GC}:N_{AT}$ results in a decrease of the spreading, as evidenced by the lower standard deviation in Fig.~\ref{rats}(b) and the higher maximum probability in Fig.~\ref{rats}(c). A linear relationship between $N$ and the mean $\alpha_0$ is observed for all ratios, with the slope of the line influenced by the ratio.
The slope $m$ for each case is: $m=0.25$ for ratio $6:1$, $m=0.45$ for $2:1$ and  $m=0.5$ for $1:1$.

\section{Conclusions}
\label{sec:con}

Motivated by the possibility that the number $\alpha$ of base pair alternations in a {circular or periodic} DNA chain might affect the dynamics of the system, we have found a {probability distribution} for this number. {Algorithms for such distributions are known for linear DNA sequences with fixed boundary conditions \cite{Robin_Daudin_1999_JApplProb_36_179}. The introduction of the periodic boundary conditions we consider in our study makes the counting of alternations a much more complicated problem due to the appearance of additional rotational and reflectional symmetries. To account for the additional complexity arising from these symmetries we have implemented P\'olya counting theory.} In particular, extending P\'olya's Enumeration Theorem for a partition-preserving group action on a partitioned set, we have constructed a well defined algorithm for calculating the number of DNA chains having a given number of alternations for particular values of the number of AT ($N_{AT}$) and GC ($N_{GC})$  base pairs.

The obtained theoretical results were compared with numerically constructed pdfs through MC simulations. We found that, in general, creating a few thousands of random DNA chains (around 5000) by MC simulations we can approximate quite accurately the theoretical pdf of $\alpha$. This means that a statistical analysis of these DNA chains will suffice to uncover the potential influence of heterogeneity on the dynamic behavior of the considered DNA model.

In addition, approximating the obtained pdfs by Gaussians we investigated the effect of the number of the two base pairs, as well as their ratio on various characteristics of the pdfs, like their mean value, their standard deviation and their maximum.

\section*{APPENDIX}
Here we present a Python computer code implementing  the algorithm of Sect.~\ref{sec:algorithm}. The function  \texttt{necklace\_count(n, B, W)} returns the total number of possible necklaces under the symmetry constraints with $2n$ alternations, $B$ black beads and $W$ white beads.

\begin{lstlisting}[language=Python]
from math import gcd
# Compute binomial coefficients in linear time.
def binomial(n, k):
    if k > n or k < 0:
        return 0
    if k = = 0:
        return 1
    if k > n//2:
        return binomial(n, n-k)
    return (n * binomial(n-1, k-1)) // k

# Compute the Euler totient function \phi(n), which
# gives the number of integers 0 < d <= n that are
# relatively prime to n.
def totient(n):
    count = 0
    for d in range(1, n+1):
        if gcd(d, n) = = 1:
            count += 1
    return count

# Get the x^r coefficient of our weight generating functions f(x^m)^n,
# where:
#    f(x) = x + x^2 + x^3 + ...
def weight_gf(r, m, n):
    if n = = 0:
        if r = = 0:
            return 1
        return 0
    if r%m != 0:
        return 0
    if (r//m) < n:
        return 0
    return binomial((r // m)-1, n-1)

# Get the x^r coefficient of a binary product of weight generating 
# functions f(x^m1)^n1 * f(x^m2)^n2, where:
#    f(x) = x + x^2 + x^3 + ...
def binary_weight_gf(r, m1, n1, m2, n2):
    total = 0
    for i in range(1, r):
        total += weight_gf(i, m1, n1) * weight_gf(r-i, m2, n2)
    return total

# Compute the number of necklaces up to dihedral symmetry with
# 2n alternations, B black beads and W white beads.
def necklace_count(n, B, W):
    # First we count the contributions from the cyclic part 
    # of the cycle index.
    count = 0
    for d in range(1, n+1):
        if n%d != 0:
            continue
        count += totient(d) * weight_gf(B, d, n//d) * 
                 weight_gf(W, d, n//d)
    # Next we count the contributions from the dihedral part 
    # of the cycle index.
    if n%2 == 0:
        count += (weight_gf(B, 2, n//2) *
                  binary_weight_gf(W, 1, 2, 2, (n-2)//2) * (n//2))
        count += (weight_gf(W, 2, n//2) *
                  binary_weight_gf(B, 1, 2, 2, (n-2)//2) * (n//2))
    else:
        count += (binary_weight_gf(B, 1, 1, 2, (n-1)//2) *
                  binary_weight_gf(W, 1, 1, 2, (n-1)//2) * n)
    return count // (2*n)
\end{lstlisting}
\section*{Acknowledgements}

M.H.~and G.P-J.~acknowledge  financial assistance from the National Research Foundation (NRF) of South Africa towards this research. G.K.~and Ch.S.~were supported by the Erasmus+/ International Credit Mobility KA107 program. Ch.S.~acknowledges support by the NRF of South Africa (IFRR and CPRR Programmes), the UCT (URC Conference Travel Grant) and thanks Hans-Peter Kunzi for useful discussions.

\end{document}